\documentclass[10pt,amssymb,twocolumn,notitlepage]{article}

\usepackage{epigraph}
\usepackage[affil-it]{authblk}
\usepackage{dcolumn}
\usepackage[stable]{footmisc}
\usepackage{blkarray}
\usepackage{bm}
\usepackage[mathscr]{euscript}
\usepackage[font=scriptsize]{caption}
\usepackage{subcaption,graphicx}
\usepackage[table]{xcolor}
\usepackage[margin=1.7cm]{geometry}
\usepackage{tcolorbox}
\usepackage{hyperref}
\usepackage{gensymb}
\usepackage[font=footnotesize]{caption}
\usepackage{hyperref}
\hypersetup{
    colorlinks=true,
    linkcolor=blue,
    filecolor=magenta,      
    urlcolor=blue,
}
\usepackage{steinmetz}
\usepackage{amssymb}
\usepackage{makecell}
\usepackage{algorithm,lipsum,xcolor,caption}
\usepackage{enumitem}
\usepackage{amsmath,esint}
\usepackage{fancyvrb}
\usepackage{xcolor}
\usepackage{arydshln}
\usepackage{mathtools}

\definecolor{newblue}{rgb}{0.0, 0.28, 0.67}
\definecolor{newgreen}{rgb}{0.13, 0.55, 0.13}
\definecolor{newred}{rgb}{0.87, 0.72, 0.53}

\newcommand{\R}{\mathbb{R}}
\newcommand{\lbar}{\{\kern-0.5ex|}
\newcommand{\rbar}{|\kern-0.5ex\}}
\usepackage{stackengine} 
\newcommand\oast{\stackMath\mathbin{\stackinset{c}{0ex}{c}{0ex}{\ast}{\bigcirc}}}

\definecolor{newblue}{rgb}{0.0, 0.28, 0.67}
\definecolor{newgreen}{rgb}{0.13, 0.55, 0.13}
\definecolor{newred}{rgb}{0.87, 0.72, 0.53}

\title{Multisets  }
\author{Luciano da Fontoura Costa \\ \emph{luciano@ifsc.usp.br}}
\affil{S\~ao Carlos Institute of Physics -- DFCM/USP} 
\date{15th Oct 2021}

\begin{document}

\twocolumn[
\begin{@twocolumnfalse}
    \maketitle
    \begin{abstract}
    Multisets are sets that allow repetition of elements, therefore accounting for their frequency
    of observation.  As such, multisets pave the way to a number of interesting possibilities
    of both theoretical and applied nature.  In the present work, after revising the main
    aspects of traditional sets, we introduce some of the main
    concepts and characteristics of multisets, which is followed by their generalization to
    take into account vectors matrices.  An 
     approach is also proposed in which the real, negative multiplicities are allowed, implying
    the multiset universe to become finite and well-defined, corresponding to the multiset
    with all support element associated with null multiplicities.  It then becomes possible 
    to define the complement operation
    in multisets in a robust manner, which allows properties involving complement  --
    including the De Morgan theorem -- to be recovered in multisets. In addition, it becomes
    possible to extend multisets to functions (which become multifunctions), scalar fields and 
    other continuous mathematical
    structure, therefore achieving an enhanced space endowed with all algebraic operations
    plus set theoretical operations including union, intersection, and complementation.   The possibility
    to define a set operation between mfunctions, namely the common product, 
    that is analogous to the traditional inner product
    is also proposed, paving the way to obtaining respective mfunction transformations, and it
    is argued that the Walsh functions provide an orthogonal basis for the mfunctions space under
    the common product.  This result also allowed the proposal of performing integrated signal processing
    operations on mset mfunctions, including filtering and enhanced template matching.
    Relationships between the cosine similarity index and the
    Jaccard index are also identified, including the presentation of an intersection-based
    variation of the cosine index.   The potential
    of multisets in pattern recognition and deep learning is also briefly characterized and
    illustrated, more specifically regarding the frequent issue of comparing clusters or
    densities.
    \end{abstract}
\end{@twocolumnfalse} \bigskip
]

\setlength{\epigraphwidth}{.49\textwidth}
\epigraph{`In the bag, seashells gathered long ago resound.'}
{\emph{LdaFC}}

\section{Introduction}

Multisets can be informally understood as sets capable of incorporation of repeated entries
of the same element (e.g.~\cite{Hein,Knuth,Blizard,Blizard2,Thangavelu,Singh}).
In a sense, they are at least as compatible with human experience than sets.
For instance it is more relevant
to know that our bag contains 4 apples than knowing that there are only apples.

In the present work, we aim at providing an introduction to this interesting area, while
also briefly covering the Jaccard index adapted to multiset, and applications to
pattern recognition.  Some new results are also described, including the possibility to
allow real, negative multiplicities, allowing the multiset universe to be composed of null
multiplicities for all support elements.  This allows the complement operation to become
stable, recovering analogous properties to all counterparts in traditional sets,
including the De Morgan theorem.

We start by reviewing the main concepts and properties of traditional sets, and then
present the concept of multisets, as well as some of their simpler properties, also
including several examples.
The challenges implied by the definition of a universe set for multisets is briefly
characterized and discussed.  It is also argued that the operations of sum and subtraction 
between multisets correspond to one of the main distinction between multiset and set theories.  

The possibility to generalize multisets to several other mathematical structures including
vectors, matrices, functions, scalar and vector fields, as well as probability densities are
approached next, including several examples.  In particular, the extension of multisets as
representations of functions and scalar fields paves the way for obtaining hybrid
expressions involving combinations of the the set operations of union and intersection
with algebraic expressions involving sum, subtraction, product and division of sets.

The interesting possibility to obtain an analogue of the traditional inner product functional,
the common product,
as well as respectively based transformations, is addressed next, with the introduction of
the concepts of mproducts and common products, which leads to the understanding of the
Jaccard index for functions as a normalized version of the common product.
The relationship between the cosine similarity index and the Jaccard index is also
addressed from the perspective of multisets and multifunctions.

We then presents how the definition of the common product between two msets or mfunctions 
paves the way to obtaining integrated signal operations including filtering and
enhanced template matching.

The Jaccard index, as well as its extension to multisets and multiple arguments, is then
briefly presented as an interesting manner to compare any of the mathematical structures
mentioned above after they have been transformed into respective multisets.

The application of multisets and the Jaccard index to quantify the relationship between
two or more clusters is then described with respect to an example related to the iris
dataset.  The measurement of the separation between clusters corresponds to an
important issue in both pattern recognition (e.g.~\cite{DudaHart,Koutrombas}), 
deep learning (e.g.~\cite{Hinton,Schmidhuber,CostaDeep}), and 
modeling (e.g.~\cite{CostaModeling,CostaAmple})

For simplicity's sake, the term multisets are henceforth abbreviated as \emph{msets}.

\section{Traditional Sets}

A \emph{set} is an unordered collection of items, or \emph{elements}, which are not
allowed to repeat.   A set $A$ with elements a, b, and c is typically represented as:
\begin{equation}
  A = \left\{ a, b, c \right\}  \nonumber
\end{equation}

The two essential properties of sets therefore are that the elements may appear in any
order, which distinguish sets from vectors, and that the elements cannot be repeated.

The number of elements in a set is called its \emph{cardinality} or \emph{size},
being represented as $|A|$.

A \emph{subset} $B$ of a given set $A$ consists of a set so that any of its elements
are contained in $A$.  If $A$ contains $N$ elements, there will be $|A| = 2^N$ possible
subsets that can be derived from it.  The set containing all possible subsets of $A$
is called its \emph{power set} $P^A$.

An important point about sets that is sometimes overlooked regards the fact that
they always refer to a respective \emph{universe set} $\Omega$.  More specifically,
once this set is established, any possible set needs to be a subset of $\Omega$.
Observe that $\Omega$ can have any type of element, though the situation where
the elements are homogenous is of particular interest.

In case some sets are given but the universe set is not provided, it is still possible to 
estimate the respective universe set (within hypotheses and subject to incompleteness
in case new sets appear) as corresponding to the union of the supplied sets.

The universe set is of fundamental importance because the operation of \emph{complement}
of a set is defined with respect to the universe set.  More specifically, the
complement of a set $A$ consists of all elements of $\Omega$ that are not part of $A$.
The complement of a set is henceforth represented as $A^C$, being implicit that
the operation refers to a given $\Omega$.

Sets can be finite or infinite, as well as discrete or continuous.  A \emph{finite set}
is any set A so that $|A| < \infty$.  A \emph{discrete set} is characterized by having
all its elements corresponding to isolated points $p$, in the sense that each of these
points possesses a neighborhood which when united with the universe set yields only $p$.
Any continuous set is infinite, but discrete sets can be finite or infinite.

An interesting point regards the relationship between an element, let's say `a' and
the set $\left\{ a \right\}$.  These two mathematical structures are not identical because
it is possible to include an element into $\left\{ a \right\}$, but not into `a'.

The \emph{empty set}, represented as $\phi = \left\{ \right\}$ is a subset of any possible set.

Given two sets $A$ and $B$, their \emph{union} consists of a third set $C$ containing all
elements from $A$ and $B$.  The \emph{intersect} of these two sets corresponds to a
set $C$ containing all elements that are in both $A$ and $B$.   A subset $B$ of $A$
can therefore be understood as to be so that $A \cap B = A$.  Any set is a subset of
itself.

The \emph{difference} between two sets $A$ and $B$, indicated as $A-B$, 
corresponds to the set $C$ containing all elements that are in $A$ but are not in $B$.

Given three sets $A$, $B$, and $C$ derived from a given $\Omega$, the following
properties are satisfied:
\begin{eqnarray} 
  A \cup A^C = \Omega  \label{eq:propset1} \\
  A \cap A^C = A \\
  A \cup \phi = A \\
  A \cap \phi = \phi \\
  A \cup A = A \\
  A \cap A = A \\
   A \cup B = B \cup A \\
   A \cap B = B \cap A \\
   A \cup (B \cup C) = (A \cup B) \cup C \\
   A \cap (B \cup C) = (A \cap B) \cup C \\
   A \cap (B \cup C) = (A \cap B) \cup (A \cap C) \\
   A \cup (B \cap C) = (A \cup B) \cap (A \cup C) \\
   (A \cup B)^C =  A^C \cap B^C   \text{ (De Morgan)}\\  
   (A \cap B)^C =  A^C \cup B^C \text{ (De Morgan)}   \label{eq:propsetn}
\end{eqnarray}

\section{Multisets} \label{sec:msets}

Basically, \emph{msets} are sets allowing the \emph{repetition} of elements, 
which is understood as  their \emph{multiplicity} or \emph{frequency}.  As with sets, the order of the
elements is immaterial.  Examples of mset include:
\begin{eqnarray}
   A =  \lbar a, a, b, b, b, d \rbar; \nonumber \\
   B = \lbar 1, 2, 1, 2, 1, 2, 1 \rbar \nonumber = \lbar 1, 1, 1, 1, 2, 2, 2 \rbar; \nonumber \\
   C = \lbar 1, a, 2, b, b, 3, c, c, c, 1, d, 2, a, a \rbar  = \nonumber \\
   =  \lbar 1, 1, 2, 2, 3, a, a, a, b, b, c, c, c \rbar; \nonumber \\
   D = \lbar a, a, b, d \rbar. \nonumber
\end{eqnarray}

Observe the different symbol adopted henceforth in this work in order to emphasize the distinction
between a traditional set ($\left\{ \right\}$), and a mset ($\lbar \rbar$).

A more compact representation of a mset $A$ can be obtained by using 2-tuple or pairs
$[a, m(a)]$, where `a' is an element and $m(a)$ it its multiplicity, i.e.~the number
of times it appear in $A$.  In the case of the above examples, we have:
\begin{eqnarray} 
   A =  \lbar a, a, b, b, b, d \rbar  =  \lbar[a,2]; [b,2]; [d,1] \rbar; \nonumber \\
   B = \lbar 1, 1, 1, 1, 2, 2, 2 \rbar = \lbar [1,4]; [2,3] \rbar; \nonumber \\
   C = \lbar 1, 1, 2, 2, 3, a, a, a, b, b, c, c, c \rbar = \nonumber \\
    = \lbar [1,2]; [2,2]; [3,1]; [a,3]; [b,2]; [c,3] \rbar; \nonumber \\
   D = \lbar a, a, b, d \rbar = \lbar [a,2]; [b,1]; [d, 1] \rbar. \label{eq:exs}
\end{eqnarray}

Though this type of representation of msets actually corresponds to a set,
because it is impossible to have two identical entries, we shall maintain the
`$\lbar \rbar$' notation in order to emphasize that a mset is being meant.

When referring to the multiplicity of an element, it is important to specify to which
mset this is being referred.  This can be done by writing $m_A(a)$, meaning
the multiplicity of the element $a$ in the mset $A$.

The property analogous to inclusion in sets can be stated as follows.  A mset $A$ 
is included in another mset $B$ whenever $m_A(a) \leq m_B(a)$.  
For instance, in the case of the examples above, we have $m_A(a) = 2$ and
$m_C(a) = 2$.

As with sets, it is particularly important to specify the universe of a mset. This can be
done in an analogous manner as with sets. As an example, let's obtain a possible
universe set for the two msets $A$ and $B$ above:
\begin{equation}
   \Omega = \left\{ a, b, d, 1, 2 \right\}
\end{equation}

It should be kept in mind that this universe is not analogous to the counterpart
in sets, as it does not actually account for the possible multiplicity of the involved
elements, which is unbound through the sum operation.

It is now possible to rewrite those two sets in a more complete, though redundant
manner, as follows:
\begin{eqnarray}
   P =   \lbar[a,2]; [b,2]; [d,1]; [1, 0]; [2,0] \rbar; \nonumber \\
   Q = \lbar [a,0]; [b,0]; [c,0]; [1,4]; [2,3] \rbar; \nonumber 
\end{eqnarray}

The \emph{support} of a given mset $A$ is defined as:
\begin{equation}
   \mathcal{S_A} = \left\{ x | x \in \Omega, m(x) > 0 \right\}
\end{equation}

For instance, the supports of the msets in Equation~\ref{eq:exs}
\begin{eqnarray}
   S_A = {a, b} \nonumber \\
   S_B = {1, 2} \nonumber \\
   S_C = {1, 2, a, b, c} \nonumber \\
   S_D = {a, b, d} \nonumber \\
\end{eqnarray}

As such, this set can be understood as containing all distinct elements
in $A$. Observe that the support set provides a useful index for 
identifying the possible elements in the respective msets.

\section{Multiset Operations}

A set $A$ is said to be \emph{included} into another set whenever:
\begin{equation}
   m_A(x) \leq m_b(x), \forall x \in A
\end{equation}

For simplicity's sake, we will indicate this operation using the same symbol
as for sets, i.e.~$A \subset B$, as the type of operation can be inferred from
$A$ and $B$ being sets or msets.

In the case of the mset examples above, we can write that $D \subset A$.

The \emph{union} $C$ of two msets $A$ and $B$ can be defined as:
\begin{eqnarray}
  C = A \cup B = \lbar [x,m_C(x)], x\in A \text{ or } x \in B \rbar, \nonumber \\
   \text{ with } m_C(x) = \max \left\{ m_A(x), m_B(x) \right\} 
\end{eqnarray}

Examples considering the msets in the beginning of Section~\ref{sec:msets}
include:
\begin{eqnarray}
  A \cup B = \lbar a, a, b, b, b, d, 1, 1, 1, 1, 2, 2, 2 \rbar = \nonumber \\
     \lbar[1,2]; [b, 2]; d,1]; 1,4]; [2,3] \rbar \nonumber \\
  A \cup D = \lbar a, a, b, d \rbar =  \lbar [a,4]; [b,4]; [d,2] \rbar \nonumber 
\end{eqnarray}

It is interesting to observe that the resulting multiplicity of each element
does \emph{not} correspond to the \emph{sum} of the respective multiplicities,
but to the \emph{maximum} between them.  This is a particularly important point
that deserves further contemplation, so we will be back to it after
presenting the concept of \emph{sum} of two msets.

Let $A$ and $B$ be msets.  The \emph{sum} of these two sets, henceforth
represented as $C = A + B$, is defined as:
\begin{eqnarray}
  C = A + B = \lbar [x,m_C(x)], x\in A \text{ or } x \in B \rbar, \nonumber \\
   \text{ with } m_C(x) = m_A(x) + m_B(x)  
\end{eqnarray}

\begin{figure}[h!]  
\begin{center}
   \includegraphics[width=0.7\linewidth]{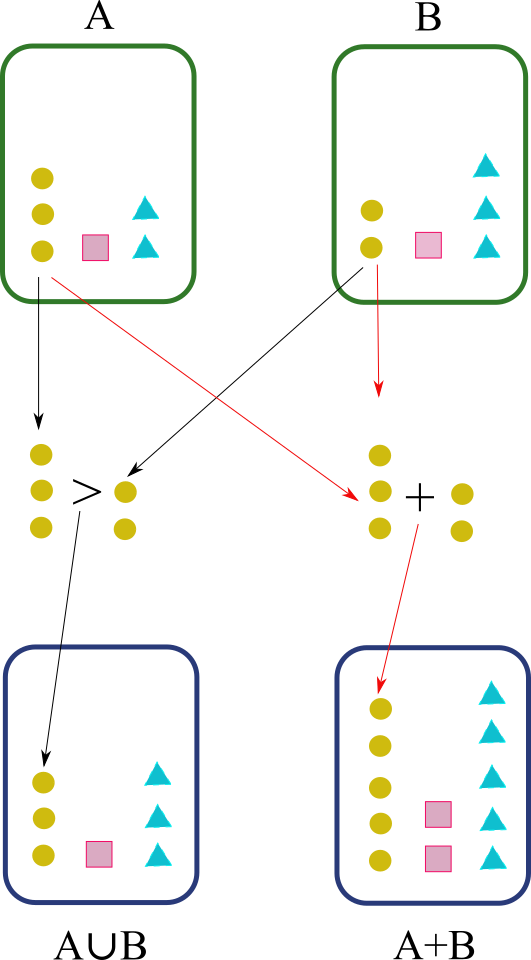}  
    \caption{The union (a) and sum (b) of two msets $A$ and $B$ typically yield quite
    different resulting msets. In the case of the union operation, each of the elements
    of the same type are compared, with the elements with the maximum multiplicity being 
    incorporated into $C$.  The sum of the two msets incorporates all the
    $m_A(x_i) + m_B(x_i)$ elements  into $C$.}
    \label{fig:uni_sum}
    \end{center}
\end{figure}

Figure~\ref{fig:uni_sum} illustrates the two different ways in which the common elements
of two msets $A$ and $B$ are collected into their respective union and sum msets.

\vspace{0.5cm}
Examples respective to the msets in the beginning of Section~\ref{sec:msets}
include:
\begin{eqnarray}
  A + B = \lbar a, a, b, b, b, d, 1, 1, 1, 1, 2, 2, 2 \rbar = \nonumber \\
     \lbar[1,2]; [b, 2]; d,1]; 1,4]; [2,3] \rbar  \nonumber \\
  A + D = \lbar a, a, a, a, b, b, b, b, d, d \rbar =  \lbar [a,4]; [b,4]; [d,2] \rbar \nonumber 
\end{eqnarray}

Thus, we have that the mset operations of union and sum are related in the
sense that both collect the elements from the two msets, but the way in which
this is done is quite different, with the multiplicities of the mset obtained by
union becoming necessarily smaller or equal than that of the mset obtained by sum,
i.e.~ $m_{A\cup B}(x_i) \leq m_{A + B}$.  

It is interesting to consider these two
operations in the context of possible respective applications.  The sum of the two
msets ensures conservation of the total number of elements (such as in 
conservative or flow-related
problems), being therefore
more indicated for related situations.  The union of two msets can be conceptually
understood as a kind of mid point between the sum  of msets and the conventional union 
of traditional sets.  

Though the union of msets will typically yield larger msets than
the union, it will not guarantee conservation of the total number of elements.
A typical situation in which the union of msets can be applied is when the incorporation
of the elements from the two msets is performed in terms of a choice, with the larger
set being taken.

The \emph{intersection} between two msets $A$ and $B$ can be defined as:
\begin{eqnarray}
  C = A \cap B = \lbar [x,m_C(x)], x\in A \text{ or } x \in B \rbar, \nonumber \\
   \text{ with } m_C(x) = \min \left\{ m_A(x), m_B(x) \right\} 
\end{eqnarray}

Examples drawn from the msets in the beginning of Section~\ref{sec:msets}
include:
\begin{eqnarray}
  A \cap B = \lbar  \rbar = \nonumber \\
  A \cap D = \lbar a, a, b \rbar =  \lbar [a,4]; [b,4]; [d,2] \rbar \nonumber \\
\end{eqnarray}

The \emph{difference} or \emph{subtraction} between two msets is expressed as:
\begin{eqnarray}
  C = A - B = \lbar [x,m_C(x)], x\in A \text{ or } x \in B \rbar, \nonumber \\
   \text{ with } m_C(x) = \max \left\{ m_A(x) - m_B(x), 0 \right\}
\end{eqnarray}

It is interesting to observe that the restriction of not having negative values can be
overlooked without great impact on the other properties and operations as addressed
in the present work.  Actually, in Section~\ref{sec:func}, we will show that allowing
negative multiplicities paves the way to a robust definition of the universe mset as
well as closed subtraction operations.

\begin{figure}[h!]  
\begin{center}
   \includegraphics[width=0.7\linewidth]{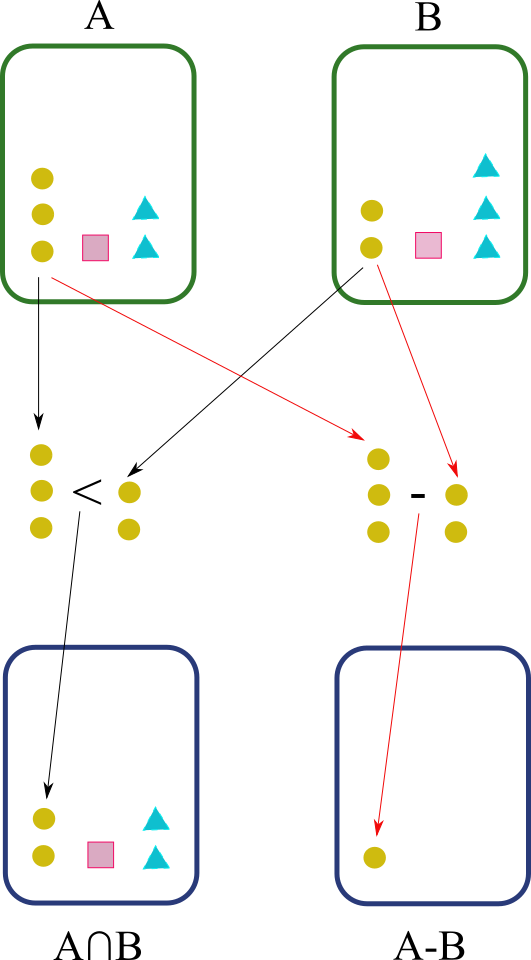}  
    \caption{The intersection (a) and difference (b) of two msets $A$ and $B$ typically yield quite
    different resulting msets. In the case of the intersection operation, each of the elements
    of the same type are compared, with the elements with the minimum respective multiplicity being 
    incorporated into $C$.  The difference between $A$ and $B$ depends on
    $m_A(x_i) - m_B(x_i)$. As the result is negative in the case of the present example,
    no elements are incorporated into $A-B$.}
    \label{fig:int_dif}
    \end{center}
\end{figure}

Figure~\ref{fig:int_dif} illustrates the intersection and difference between two msets $A$ and $B$.

\vspace{0.5cm}

Respective examples include:
\begin{eqnarray}
  A - D = \lbar b,b \rbar = \lbar [b,2] \rbar  \nonumber \\
  D - A = \lbar  \rbar  \nonumber \\
\end{eqnarray}

\section{Multisets Properties}

It can be shown that msets satisfy the following properties:
\begin{eqnarray}
  A \cup \phi = A \\
  A \cap \phi = \phi \\
  A \cup A = A \\
  A \cap A = A \\
   A \cup B = B \cup A \\
   A \cap B = B \cap A \\
   A \cup (B \cup C) = (A \cup B) \cup C \\
   A \cap (B \cup C) = (A \cap B) \cup C \\
   A \cap (B \cup C) = (A \cap B) \cup (A \cap C) \\
   A \cup (B \cap C) = (A \cup B) \cap (A \cup C) 
\end{eqnarray}

Informally speaking, the above properties can ben understood by identifying each
repeated element in an mset with unique respective tags, in which case they
would behave with respect to the operation of union and intersection exactly 
in the same way as common sets.

So, we have that msets follow all the properties in Equations~\ref{eq:propset1} to Equations~\ref{eq:propsetn} , except those involving complementation.   Indeed, the definition of the complement of an mset
has been a challenging issue (e.g.~\cite{Singh}), which has to do with the fact that there is no bound
to the size of possible msets generated by making additions between any non-empty
mset.  For instance, we can write:
\begin{eqnarray}
   A = \lbar a \rbar   \nonumber \\
   A = A + A = \lbar a, a \rbar  \nonumber \\
   A = A + A = \lbar a, a, a,a \rbar  \nonumber \\
   \ldots 
\end{eqnarray}

For this reason, the useful De Morgan properties, as well as other related results,
do not hold for msets.   In addition, there are relative few properties involving the sum and
subtraction of msets.  However, when negative values are allowed for the subtraction
between msets, the universe mset consists of taking all the multiplicities as zero,
and the complement operation becomes the change of sign of the multiplicities
(see Equation~\ref{eq:muniverse}).  When applied to functions
represented as msets, this means that the null function is the respective universe.

In a sense, it is these two operations that differentiates
msets  from sets because, as commented above, msets behave like
sets respectively to their union and intersection when the elements are tagged.

\section{Multisets, Vectors, Matrices}

In this section we will discuss the relationship between msets and vectors.  
First, we recall that the elements in a vector are expected to follow a well-determined
order as indicated by their indices.  For instance, in the case of the vector in $\R^5$:
\begin{equation}
  \vec{v} = [3, 2.5, \pi , 0, -1]  \nonumber
\end{equation}

we have five indices $i = 1, 2, \ldots, 5$, so that we can specify the respective element values
as $v[1] = 3$, $v[2] = 2.5$, $v[3] = \pi $, $v[4] = 0$, and $v[5] = -1$.

By understanding the values of the components of a vector as generalized multiplicities,
we can immediately derive the following mset from the above vector:
\begin{equation}
  V = \lbar [1, 3]; [2, 2.5]; [3,\pi]; [4,0]; [5,-1] \rbar  \nonumber
\end{equation}

Therefore, we have that an mset can be derived from any vector, but that a vector can
be obtained from an mset only if their elements are ordered in some manner, e.g.~by
taking their respective values instead of understanding them as labels.    
This situation becomes more evident when one
considers non-numeric elements.
As such, msets can be used to study the elements
of vectors without taking into account their relative position along the vector.

It is also interesting to contemplate the relationship between the above discussion and
the traditional sets containing multiplicities.  More specifically, we can write the set containing all multiplicities in the 
vector $\vec{v}$ above as:
\begin{equation}
  \tilde{V} = \left\{ 3,  2.5, \pi , 0 ,-1\right\}  = \left\{ 0, 3,  2.5, -1, \pi \right\}  = \text{etc.} \nonumber
\end{equation}

Though $V$ and $\tilde{V}$ are very similar, they are not identical because in $V$ the
\emph{correspondence} of the elements and the respective multiplicity is maintained.  This difference
becomes critically important in case we want to apply the operations of addition and
subtraction between msets, which would be otherwise impossible in case of sets because
we would not know which element of one set should be added to which element in the
other set.

By representing vectors as msets, we not only preserve the operations of subtraction and
difference, but also incorporates the possibilities of defining intersections and unions
between any two vectors.

Another interesting possibility is to incorporate new operators for multiplication and
division into the mset framework, which can be done straightforwardly, while avoiding
divisions by zero.

Interestingly, it is also possible to obtain mset representations from matrices or even
other more sophisticated mathematical structures as tensors.  In the case of matrices,
this  can be done by mapping the indices $i = 1, 2, \ldots, N_i$ and $j = 1, 2, \ldots, N_j$ 
into a single index $k$, e.g.:
\begin{equation}
   k \;  \longleftrightarrow \; N_i (j-1) + i - 1
\end{equation}
 
so that the matrix becomes a vector, which can then be transformed into the respective 
mset as described above.   Interestingly, observe that
though the obtained msets would not directly provide the respective indices of the elements,
they can be nevertheless recovered from the unified index.

\section{Functions and Scalar Fields}  \label{sec:func}

The possibility to represent vectors as msets opens the way to a number of interesting
possibilities.  One of them is to represent discrete and continuous \emph{functions} and
\emph{scalar fields} (vector fields can be approached as vectors of scalar fields).
We develop these possibilities in the following.

\begin{figure}[h!]  
\begin{center}
   \includegraphics[width=0.7\linewidth]{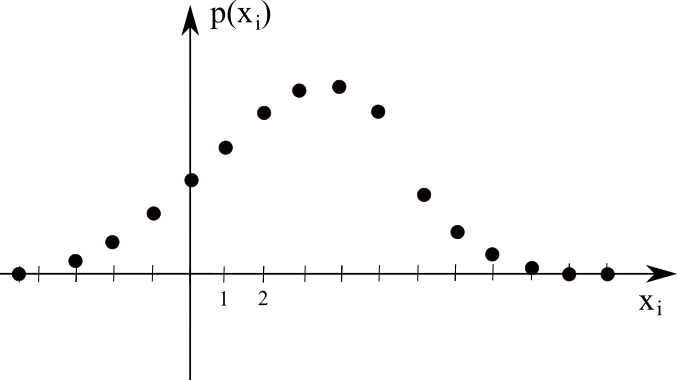}  
    \caption{A generic discrete function $f(x_i)$.}
    \label{fig:discr_func}
    \end{center}
\end{figure}

We start with a discrete function as illustrated in Figure~\ref{fig:discr_func}.

\vspace{0.5cm}

It is often interesting to represent such discrete functions in terms of sums
of Dirac delta functions.

Provided that we limit the extension of this function along the horizontal axis
(related to the function support),
e.g.~$i = 1, 2, \ldots, N$, it can be immediately represented in terms of the
vector:
\begin{equation}
  \vec{f} = [f(x_i)]
\end{equation}

which can then be transformed into the respective mset as described in the
previous section.

A similar approach can be applied to transform discrete scalar fields defined on more than
one variables into respective msets, involving the index mapping described above.

Now, we can approach the case of continuous functions and scalar fields simply by
taking the separation between the points along the horizontal axis to the limit of 0.
Observe that a homeomorphism is establsiehd between the function space and
its representation in the multisets, as provided by the bijective associations between
the function abscissa and the elements, therefore preserving the topology of the
representation.  In addition, all the operations required
for the multisets apply irrespectively of the neighborhood of the abscissa neighborhood.
Actually, it is hard to think of a bound function that will not map into the respective
multiset elements.  To any extent, the results in this work are understood to be
applicable to mset representations of functions, i.e.~mfunctions.

As a consequence, the functions and scalar fields will become associated to msets
with infinite elements, but these can still be operated by the `$\max()$', `$\min()$',
`$+$' and `$-$ operations required by the mset operations, as well as any other
function or transformation applicable to functions.  For simplicity's sake, we
shall call the multisets derived from functions as \emph{mfunctions}.

Observe, though, that all information about a function is preserved into the respective mfunction,
allowing its reconstruction.

Let's illustrate the above concepts and possibilities in terms of the two following functions:
\begin{eqnarray} \label{eq:gh}
    g(x) =  e^{-10 x^2}  \nonumber \\
    h(x) =   2 e^{-10 |x-0.1|} 
\end{eqnarray}

Figure~\ref{fig:operations} illustrates the two above functions as well as their
union, intersection, sum, and subtraction after having been converted into
respective msets.

\begin{figure}[h!]  
\begin{center}
   \includegraphics[width=1\linewidth]{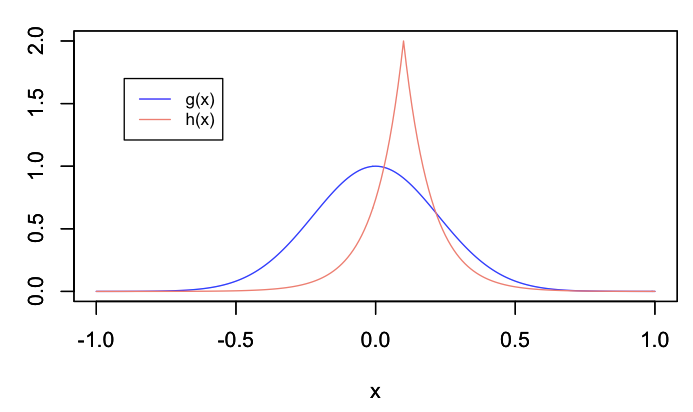}  \\ (a) \\ \vspace{0.3cm}
   \includegraphics[width=1\linewidth]{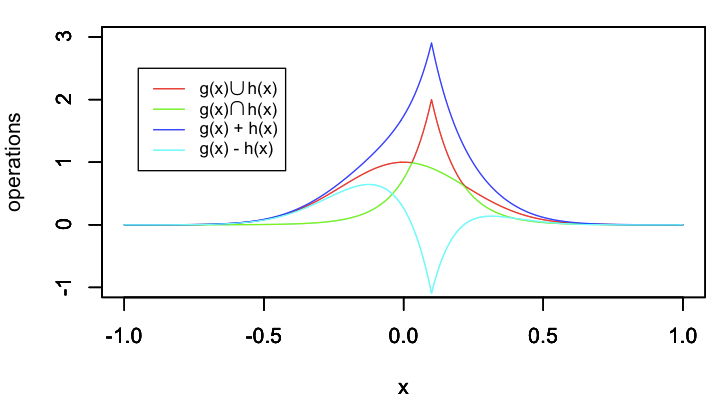}  \\  (b) \\
    \caption{Two continuous functions $g(x)$ and $h(x)$ of a single variable (a), and their respective mset operations
    (b) of union, intersection, sum, and subtraction.}
    \label{fig:operations}
    \end{center}
\end{figure}
\vspace{0.5cm}

This examples illustrate several interesting points.  First, we have the transformation of functions
into respective msets immediately allow them to be operated by union and intersection.
In addition, we observe an example that the sum of two functions is larger or equal than their union,
as well as the possibility of the subtraction operation yielding negative values.

Figure~\ref{fig:operations} shows the two additional operations of product and quotient between
the two functions $g(x)$ and $h(x)$ above, avoiding the divisions by zero.

\begin{figure}[h!]  
\begin{center}
   \includegraphics[width=1\linewidth]{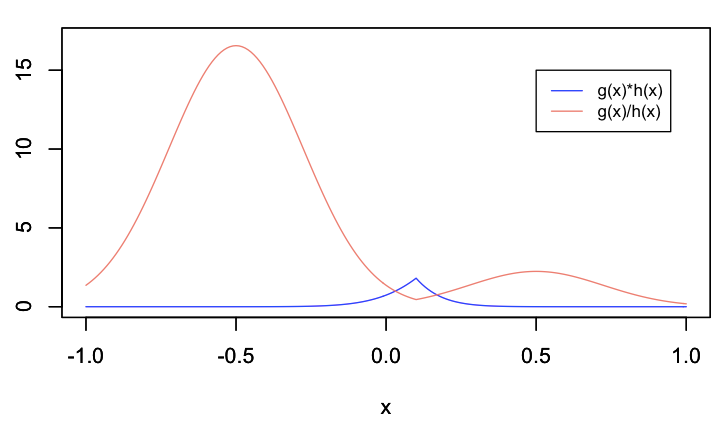}  
    \caption{Two functions $g(x)$ and $h(x)$ of a single variable (a), and their respective mset operations
    (b) of union, intersection, sum, and subtraction.}
    \label{fig:operations}
    \end{center}
\end{figure}
\vspace{0.5cm}

It is interesting to observe that the potential of the operations of union and intersection in
producing sharp derivatives and discontinuities, which contribute an interesting manner of
representing an ample range of function types as combinations of these operations, not to
mention the operations of sum, subtraction, product and quotient.

Consequently, it becomes an interesting possibility to develop transformations of functions,
analogous to the Fourier transform, considering not only series of basis functions, but also
intersections and/or unions, and/or other possible hybrid operations between msets.  One
particularly interesting benefit would be to become able to express functions with
discontinuous derivatives as combinations of functions that are completely smooth.
Also, it should be observed that the operations of sum and subtraction are bilinear,
while the minimum and maximum are not.

Indeed, the above developments also allow new functions to be obtained as combinations of logic
operations as the union and intersection and numeric operations as sum, subtraction, product,
and division.   For instance, it becomes possible to write things such as:
\begin{eqnarray}
    r(x) = (g(x)  \cap  h(x))  +   g(x) \nonumber \\
    s(x) = (g(x)  +  h(x))  \cup ( g(x)  -  h(x) )\nonumber \\
    t(x) = [g(x)  \cap   h(x)]  -  [g(x)   \cup  h(x)]  \nonumber 
\end{eqnarray}

\begin{figure}[h!]  
\begin{center}
   \includegraphics[width=0.8\linewidth]{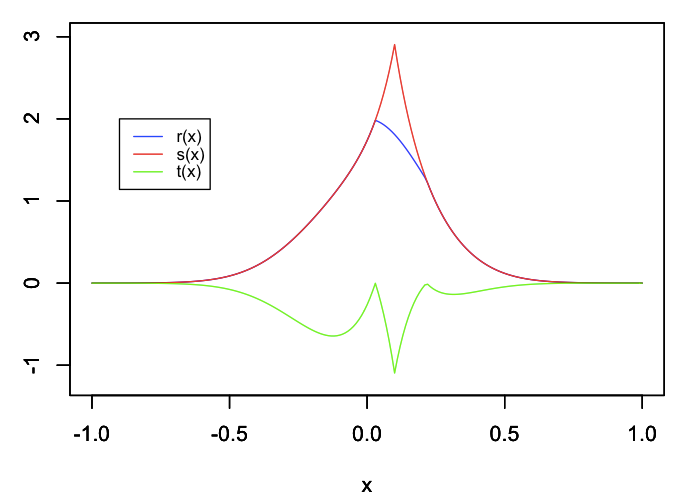}  
    \caption{The functions $r(x)$, $s(x)$ and $t(x)$ obtained through mset operations.}
    \label{fig:three_funcs}
    \end{center}
\end{figure}

These three functions are illustrated in Figure~\ref{fig:three_funcs} assuming the function in
Equation~\ref{eq:gh}.

\vspace{0.5cm}

Seen from this perspective, the complement of a multset becomes the operation of
sign change, in the sense that $r(x)^C$ would become $-r(x)$.  Indeed, it can be
verified that the De Morgan properties hold in this case, i.e:
\begin{eqnarray}
   - [g(x) \cap h(x)] = -g(x) \cup -h(x) \nonumber \\
   - [g(x) \cup h(x)] = -g(x) \cap -h(x) \nonumber 
\end{eqnarray}

The potential of these hybrid constructions is as large as our imagination, because it allows
the incorporation of much of the concepts, structures, and properties of both arithmetic,
set theory, and also logic (which is directly analogous to set theory).

It should also be observed that the above results could 
only be  achieved from allowing negative values for the multiplicity in the
difference operation between msets.  The sign change therefore recovers many properties
analogous to those involving complements of sets, and more.  As a consequence, the universe mset 
becomes as follows:
\begin{equation}
  \Omega = \lbar [x_1,0]; [x_2,0]; \ldots; [x_N,0] \rbar  \label{eq:muniverse}
\end{equation}

where $x_i$ are the elements of the support of the msets.

Other properties analogous to those of traditional sets are also consequently recovered, including
$A \cup A^C$.

\section{The Multiset Jaccard Indices}

The Jaccard index represents an effective and conceptually appealing manner to quantify
the similarity between any two sets $A$ and $B$ (e.g.~\cite{Jaccard, Jac:wiki, CostaJaccard}), 
having therefore being extensively applied in a vast range of problems in several scientific and technological
fields.

In its most basic form, the Jaccard index between an two sets $A$ and $B$ can be
expressed as:
\begin{equation}
   \mathcal{J}(A,B) = \frac{|A \cap B|}   {|A \cup B|}
\end{equation}

It is possible to adapt the Jaccard index to msets by making:
\begin{equation}  \label{eq:Jac_multi}
   \mathcal{J}_M(A,B) =   \frac{\sum_{i=1}^N\min{(m(a_i),m(b_i))}}  {\sum_{i=1}^N\max{(m(a_i),m(b_i))}} 
 \end{equation}
 
where $a_i$ and $b_i$ are the elements of the sets $A$ and $B$, respectively, and N
is the cardinality of the universe of those sets.  We also have that $0 \leq \mathcal{J}_M(A,B) \leq 1$.
 
As an example, let's consider $A = \left\{ a, b, b, c, c, c \right\}$ and $B = \left\{ a, b, c, c, d \right\}$.
Then, we have:
\begin{equation}  
   \mathcal{J}(A,B) =   \frac{1 + 1 + 2 + 0}  {2 + 3 + 5 + 1} = \frac{4} {11} 
 \end{equation}
  
 It is possible to adapt the Jaccard index to mfunctions by making:
\begin{equation}  \label{eq:Jac_multi}
   \mathcal{J}(A,B) =   \frac{\int_{\Phi}\min{(m_A(\vec{x}),m_B(\vec{x}))}}  {\int_{\Phi}\max{(m_A(\vec{x}),m_B(\vec{x}))}} 
 \end{equation}

where $\Phi$ is the common support of the two functions or scalar fields, and $0 \leq \mathcal{J}(A,B) \leq 1$.
 
As such, the Jaccard index can be understood as a \emph{functional}, or \emph{mfunctional}  
of the two functions of scalar fields. 

The Jaccard index has been enhanced and extended to functions, scalar fields, joint variations and more than
2 sets~\cite{CostaJaccard}.  In particular, the latter type of Jaccard index for 3 sets $A$, ,$B$ and
$C$ can be written as:
\begin{equation}  \label{eq:Jac3}
   \mathcal{J}(A,B,C) = \frac{A \cap B \cap C}   {A \cup B \cup C}
\end{equation}

\section{Inner Product and Transforms}  \label{sec:inner}

The concept of multifunctions as proposed above has some interesting implications.
In particular, it effectively endows mfunctions and scalar fields with the analogous of the
set operations of union and intersection, which immediately motivates several important
questions such as what are the properties involving both set and algebraic properties of
mfunctions.  Another possibility of particular relevance because of its immediate
relationship with the concepts of similarity and transformation regards if it is possible
to define a set-related operation between two mfunctions which would be analogous to
the \emph{inner product}.   In this section we propose one such operation and then
use it to define transforms of mfunctions in a manner similar to Fourier transform, but
without the property of orthogonality extending over all distinct sinusoidals.

We start by recalling the inner product between two real mfunctions $f(x)$ and $g(x)$:
\begin{equation}
   \left< f(x), g(x) \right> = \int_{-\infty}^{\infty} f(x) g(x) dx
\end{equation}

Any two functions $f(x)$ and $g(x)$ are said to be \emph{orthogonal} if and only if
their inner product is zero.  The prototypical example of orthogonal functions consists
in sinusoidals $s(x, \omega, \phi) = sin(\omega x + \phi)$, with $\omega= 2 \pi f $.  It is known that:
\begin{equation}
   \left< f(x), g(x) \right> = \int_{-\infty}^{\infty} s(x, \omega, \phi) s(x, \omega, \tilde{\phi})  dx = 0
\end{equation}

provided $\phi \neq \tilde{\phi}$.    However, when $\phi \neq \tilde{\phi}$, we 
have $\left< f(x), g(x) \right>$=1.  Therefore sinusoidals are orthogonal.

Now, the similarity between a generic function $f(x)$ and a sinusoidal $s(x, \omega, \phi)$
can be readily expressed by the respective inner product, i.e.:
\begin{equation}
   \left< f(x), s(x, \omega, \phi) \right> = \int_{-\infty}^{\infty} f(x) s(x, \omega, \tilde{\phi})  dx = 0
\end{equation}

The higher the obtained value, the more similar the function $f(x)$ is to the sinusoidal.
Actually, the inner product $\left< f(x), s(x, \omega, \phi) \right>$ corresponds directly to the
Fourier transform of $g(x)$, which provides an effective manner of expressing $f(x)$ as 
a linear combination of sinusoidals, which are more effectively handled as complex
exponentials.  

A particularly important property of the Fourier, as well as other orthogonal transforms is that,
given the orthogonality between the basis functions, the calculation of the coefficients can
be made independently one another. 

Having obtained the transformation coefficients for every possible sinusoidal, it is then
possible to recover the original function $f(x)$. 

Now, we develop an analogous construction regarding the set operations incorporated into
mfunctions.   

First, we start by defining the mfunctional corresponding to the integral of the intersection
between any two functions $f(x)$ and $g(x)$.  It is proposed here that this can be done as follows:
\begin{equation}
   \ll f(x), g(x) \gg = \int_{-\infty}^{\infty}  s_f s_g \min(s_f f(x), s_g g(x)) dx 
\end{equation}

where $s_f$ and $s_g$ are the signs of $f(x)$ and $g(x)$, and $p(x)  = s_f s_g \min(s_f f, s_g, g)$
is henceforth called the \emph{mproduct} of $f(x)$ and $g(x)$.  In order to distinguish between
the traditional and set-based inner products, we will refer to the latter as \emph{common product}.

This integral can be verified to effectively correspond to the common area between
the two functions, with the area being taken considering the functions signs as in the above expression.

An interesting point here concerns the fact that the above integral requires the elements of the
mfunctions to be taken in a pre-specified order, which does not need to correspond to the
real line.  As such, it becomes possible to apply the common product to quantify the similarity
between any two msets in which some order relationship can be established along
their shared support.  For instance, the elements could correspond to the alphabet letters,
dates of historic events, and so on.

It is also possible to define other functionals, such as:
\begin{equation}
    f(x) \oast  g(x)  = \int_{-\infty}^{\infty}   \max(s_f f(x), s_g, g(x)) dx 
\end{equation}

henceforth called the \emph{sup product}.

Thus, the Jaccard index for mfunctions can now be expressed as:
\begin{equation}
    \mathcal{J}(f(x) ,  g(x))  = \frac{ \ll f(x), g(x) \gg   }   {    f(x) \oast  g(x)  }
\end{equation}

As such, the Jaccard index can be understood as corresponding to the common product,
which quantifies the similarity between the two msets or mfunctions in terms of the
intersection-based common product (a functional), normalized relatively to the union
of the mfunctions.

Let's illustrate this mfunctional respectively to sinusoidals.  For simplicity's sake, we
shall consider that all functions share the same support $[0, T]$, $T = 1f$.

In case $f(x) = g(x) = sin(\omega x)$, we immediately have that their intersection is 
$f \cap g = g = g$, from which
\begin{eqnarray}
   \ll f(x), g(x) \gg =  \nonumber \\
    =\int_{0}^{T/2}  sin(\omega x) dt  + \int_{T/2}^{T}  sin(\omega x) dt  = 1  \nonumber
\end{eqnarray}

Now, in case $f(x) = - g(x) = \sin(\omega x)$, we have $f(x) \cap g(x) = 0$, implying
\begin{eqnarray}
   \ll f(x), g(x) \gg =\nonumber \\
    =- \int_{0}^{T/2}  sin(\omega x) dt  - \int_{T/2}^{T}  sin(\omega x) dt  = -1  \nonumber
\end{eqnarray}

For $f(x) = \sin(\omega)$ and $g(x) = \cos(\omega)$, it follows that:
\begin{eqnarray}
   \ll f(x), g(x) \gg =\nonumber \\
   = \int_{0}^{T/4}  sin(\omega x) dt  -  \int_{T/4}^{T/2}  sin(\omega x) dt  + \nonumber \\
       + \int_{T/2}^{3T/4}  sin(\omega x) dt  -  \int_{3T/4}^{T}  sin(\omega x) dt  = 0 
\end{eqnarray}

So far, the common product between the sinusoidals has presented a remarkable similarity
with the traditional inner product, being identical for the three situations above.  In addition,
we can conclude that the sine and cosine are orthogonal also regarding the common product.

However, let's now consider that $f(x) = \sin(\omega)$ and $g(x) = \sin(\omega + \phi)$, with $\phi\neq 0$.
Unlike the traditional inner product, the common product of these two functions will no longer
be zero, being nevertheless comprised in the interval $(0,1)$. 
As a consequence, except for the case $f(x) = \sin(\omega)$ and $g(x) = \cos(\omega)$,
no other combination of sinusoidals will result orthogonal one another.

Nevertheless, it is still possible to define transformations capable of decomposing
(and recovering) a generic mfunction $f(x)$ in terms of a
combination of basic mfunctions such as sinusoidals.  This can be readily accomplished as
follows.  

Let $g_i(x)$, $i = 1, 2, \ldots, N$, be the basis function of the transformation, and $f(x)$ be the 
function to be analyzed.    The transformation coefficients can be obtained as:
\begin{align}
   &\text{for } i = 1, 2, \ldots, N \nonumber \\
   &\quad \mathfrak{c}_i = \ll f(x), g_i(x) \gg = \nonumber \\
    &=\quad \int_{-\infty}^{\infty}  s_f s_{g_i} \min(s_f f(x), s_{g_i}, g(x)) dx;
   \nonumber \\
   &\quad f(x) = f(x) - g_i(x)  \nonumber
\end{align}

Each of these coefficients corresponds to the effectively shared area between each pair of
functions.  Also, it is important to normalize the function being transformed so that it has
area equal to one.

In the case of sinusoidal mfunctions, because they are generally not orthogonal, we will have
that the above coefficients will depend on the order in which the basis functions are applied,
so that several distinct transformations can be derived for a same function and a same basis.
However, it is possible to employ optimization methods in order to achieve specific objectives.

The original function $f(x)$ can then be recovered to a given precision (depending on the
number and type of basis functions) as:
\begin{align}
   & f(x) = 0; \nonumber \\
   &\text{for } i = 1, 2, \ldots, N \nonumber \\
   &\quad f(x) = f(x) + \mathfrak{c}_i g_i(x)
\end{align}

As a matter of fact, the Walsh functions (e.g.~\cite{Harmuth,Tzafestas}) can be verified to
provide an orthogonal basis for the multifunctions with real values under the common product
binary operation.  Because these functions also provide a basis in real vector spaces, it
represents a structural characteristic shared by both multifunctions and functions, therefore
defining an additional link between these two spaces.

\section{The Cosine Similarity from the Multiset Perspective}

Given two vectors $\vec{X}$ and $\vec{Y}$ with $N$ elements each, 
their cosine similarity  is commonly expressed as:
\begin{equation}
   \text{Cs}(\vec{X}, \vec{Y}) = Cos(\theta)  =  \frac{1}{|X| |Y|}  \sum_{i = 1}^N   X_i Y_i 
\end{equation}

where $\theta$ is the smallest angle between the two vectors.
Typically the cosine similarity assumes all the elements to be nonnegative.

Now, we have already seem that vectors and multisets are closely related,
in the sense that the latter can be derived from the former.  Consequently,
the cosine similarity index can be immediately understood in terms of multisets.

Let $X$ and $B$ be the multisets derived from $\vec{X}$ and $\vec{Y}$,
therefore sharing the same support with $N$ elements.  Then, we can write:
\begin{equation}
   \text{Cs}(X,Y)  =  \frac{1}{\sum_{i=1}^N X_i  \; \sum_{i=1}^N Y_i  } \sum_{i = 1}^N   m(X_i) m(Y_i) 
\end{equation}

Therefore, the cosine similarity index can be understood as a normalized version of the Jaccard index
where the union between the multisets has been replaced by their respective product, 
normalized by the product of the number of elements in the two multisets.

The above reasoning extends directly to two real mfunctions $f(x)$ and $g(x)$:
\begin{align}
  &\text{Cs}(f(x),g(x))  =  \nonumber \\
  &\quad =  \frac{1}{\int_{\Phi} |f(x)| dx \; \int_{\Phi} |f(x)| dx} 
     \int_{\Phi}  f(x) g(x) dx \nonumber
\end{align}

Now, it becomes possible to adapt the cosine similarity to intersection-based similarity by
adopting the mproduct proposed in Section~\ref{sec:inner}, i.e.:
\begin{align}
   &\text{Cs}(f(x),g(x))  =  \nonumber \\
   & =  \frac{1}{\int_{\Phi} |f(x)| dx \; \int_{\Phi} |f(x)| dx} 
     \int_{\Phi}  s_f s_g \min(s_f f(x), s_g, g(x)) dx \nonumber
\end{align}

\section{Toward Integrated Signal Processing}

The definition of the set common product as an analogue counterpart of the inner product
binary operation between two mfunctions paves the way to  the handling and modification
of signals by taking into account set operations including union and intersection \emph{together}
with all other algebraic operations such as sum, subtraction, product, division, etc.~of mfunctions.
The resulting ensemble of concepts and methods can then constitute what is henceforth called
\emph{integrated signal processing}.  This interesting subject is developed in the current section.

We have already seen that the set common product plays a role analogous to the inner product.
This enables us to define the \emph{mconvolution}, namely the convolution between two
multifunctions (also extended to other msets) as:
\begin{equation}
  f(x) \circ g(x) [y] = \int_{-\infty}^{\infty}  \ll f(x) g(y-x) \gg dx
\end{equation}

with $y \in \R$.    Observe that, as with the traditional function convolution, the mconvolution is commutative.

The \emph{mcorrelation} between two mfunctions can then be expressed as:
\begin{equation}
  f(x) \diamond g(x) [y] = \int_{-\infty}^{\infty}  \ll f(x) g(x-y) \gg dx
\end{equation}

In case $g(x)$ is even, the mcorrelation becomes closely related to the mconvolution.
It can be verified that the mcorrelation quantifies the intersections areas between the functions,
suitably adapted to take into account negative multiplicities.

Because the mcorrelation does not penalize the difference between the two functions at a given
position $y$ with respect to the areas that are not overlapped, it becomes interesting to define
a \emph{set similarity mcorrelation}, henceforth referred to as \emph{scorr}, as follows:
\begin{equation}
  f(x) \Box g(x) [y] = \int_{-\infty}^{\infty}  \frac{ \ll f(x) g(x-y) \gg } {f(x) \oast g(x-y)} dx
\end{equation}

which, interestingly, is precisely the same as:
\begin{equation}
  f(x) \Box g(x) [y] = \int_{-\infty}^{\infty}  \mathcal{J} (f(x),g(x-y)) dx
\end{equation}

Thus, in case a more strict match is desired, the similarity correlation should be used instead
of the mcorrelation.  
Observe that the above approach extends immediately to scalar fields of any dimension as
well as  to generalizations of the Jaccard index as those described in~\cite{CostaJaccard}.

The definition of the above mconvolutions and mcorrelations paves the way to 
theoretical and applied developments in
all areas that already employ the traditional function convolution, including filtering, control theory,
signal and image analysis, template matching, to mention but a few possibilities.  This potential
is illustrated with respect to template matching in the following.

The operation of \emph{template matching} consists of, given a function $f(x)$, to quantify,
along $X$,  the similarity between its portions and another
reference template function $g(x)$.  This can be immediately implemented by by applying the 
mcorrelation or similarity mcorrelation to those two functions.  High resulting values indicate
portions of $f(x)$ that are similar to $g(x)$.

\begin{figure*}[h!]  
\begin{center}
   \includegraphics[width=0.7\linewidth]{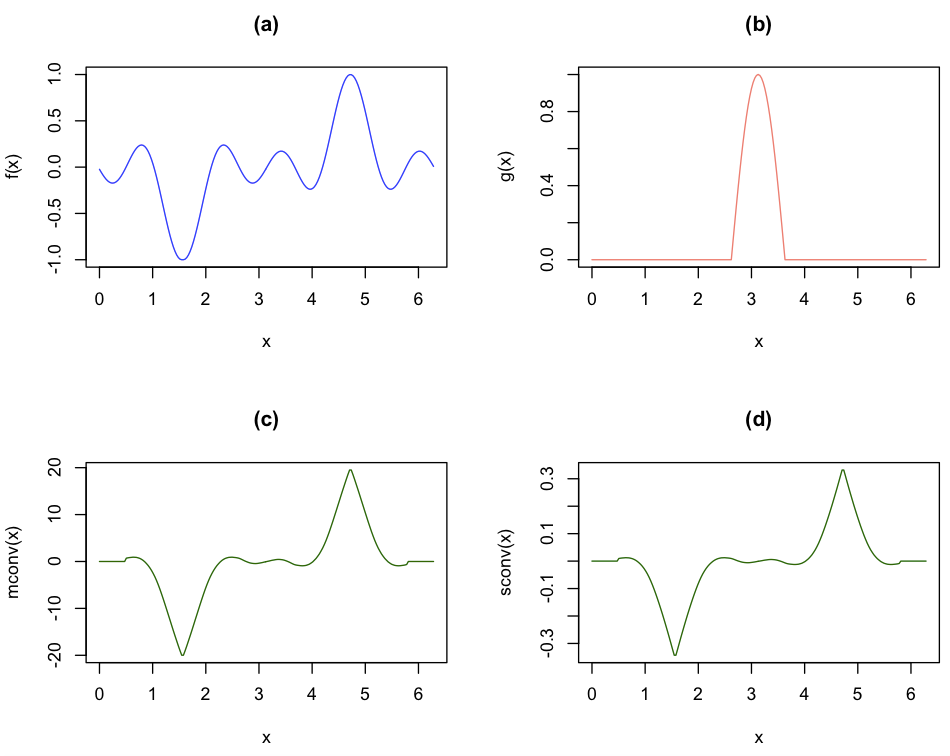}  
    \caption{Template matching through mcorrelation and scorrelation.  
    The templated in (b) is to be compared
    with the mfunction in (a).  This can ben effectively achieved by using the mcorrelation
    between these two functions, whose result is shown in (c).  The scorrelation (actually
    the Jaccard index for continuous functions) can also be employed in case a more strict quantification of the
    local similarity is required (d).}
    \label{fig:tempmatch}
    \end{center}
\end{figure*}

Figure~\ref{fig:tempmatch} presents the result of matching the template in (b) with the function
in (a) by using mcorrelation (c) and similarity mcorrelation.

\vspace{0.5cm}

As it can be verified, both the mconvolution and the sconvolution yielded a precise identification
of the maximum similarity between the portions of function $f(x)$ with the template function $g(x)$,
not only for the positive parts, but also with respect to the negative.
The secondary matchings appeared with substantially smaller values. 

\begin{figure}[h!]  
\begin{center}
   \includegraphics[width=0.7\linewidth]{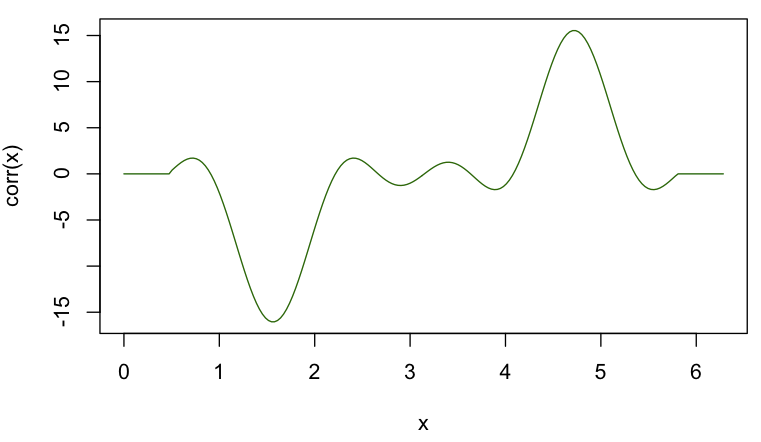}  
    \caption{The standard cross-correlation obtained for the two functionsin Figure~\ref{fig:tempmatch}.
    It can be observed that not only less sharp peaks are obtained, but also that the secondary
    peaks are substantially higher than in the results in Figure~\ref{fig:tempmatch}.}
    \label{fig:corr}
    \end{center}
\end{figure}

For comparison purposes, Figure~\ref{fig:corr} presents the result of the standard cross-correlation
function between the two functions in Figure~\ref{fig:tempmatch}.

\vspace{0.5cm}

The advantages of the mcorrelation and scorrelation relatively to the standard correlation
are evident.  Not only the main detected peaks are sharper, but also the secondary peaks
are substantially smaller in the case of the mcorrelation and scorrelation.  Actually, the
peaks obtained by the correlation result are even wider than in the original function.

For these reasons, and also considering additional experimental results by the author
(to be published), it is proposed here that the mcorrelation and scorrelation has potential
for enhanced template matching results, thanks to the sharper and
more distinguished peaks.  For similar reasons, it is also believe that
more commensurate results will be obtained, in general, in tasks such as measuring the 
similarity between  clusters and estimating joint variation of random variables 
by using the common product (see~\cite{CostaJaccard}).
The adoption of the coincidence index described in that same reference also tends
to improve the performance of the mcorrelation and scorrelation further.

Some considerations are also due regarding the computational expenses implied by the mcorrelations
and scorrelations.  One important point here is that these mfunctions operations are particularly
cheap, as they involve just the comparison between two points as implied by the minimum
and maximum operations.  Though the traditional convolution and correlation can be performed
more effectively in the Fourier domain, there will still be $log(N)$ levels and complex products
to be executed.

Given that the operation of neurons is often understood as being related to convolution or
correlation, it becomes particularly interesting to adopt the mconvolution and sconvolution for
implementing and modeling the linear portion of neurons in areas such as computational
neuroscience and deep learning.

\section{Generalized Multisets and Logic}

Set operations are naturally associated to Boolean expressions.  For instance, given
three sets $A$, $B$, and $C$,  consider:
\begin{align*}
   \textbf{Dataset Domain }  &\Longleftrightarrow  \textbf{Model Domain } \nonumber \\
   A   =  B^C    &\Longleftrightarrow  m_A = \neg m_B  \nonumber \\
   A   = B \cup C  &\Longleftrightarrow m_A = m_B  \lor m_C  \nonumber \\
   A   = B \cap C  &\Longleftrightarrow m_A = m_B  \land m_C  \nonumber \\
   A   = B - C  =  A \cap B^C   &\Longleftrightarrow m_A = m_B  \land \neg m_C  \nonumber 
\end{align*}

where $m_A$, $m_B$ and $m_C$ are models associated to the datasets $A$, $B$, and $C$.
 
It is this intrinsic association, so natural to humans to the point of sometimes not being realized,
that often provides the basis for developing new models (logical combinations of models)
while considering the adherence between respective datasets.

The generalization of the concept of multisets to several mathematical structures, as well as the
identification of their relationship with several types of similarity indices, paves the way to 
applying multisets and mfunctions effectively in modeling activities.  

As sets are naturally expanded to virtually every  other mathematical structure, it becomes
a  topic of particular interest to identify which are the logical constructs corresponding to 
multiset operations such as sum, subtraction, product and division.  Probably, these logic
concepts have to do with the incorporation of the multiplicities implied by multisets into
the logic reasoning, suggesting a logic with weighted or graded Boolean variables.  For instance,
when we say ``two apples united with three pears'' (multiset union), which has the direct logic meaning of
``apples and pears'', when transformed to ``two apples plus three pears'' (multiset sum) 
would mean ``two apples and three pears''.  

The three levels of purported relationships between two sets, propositions, or functions can
be summarized as in Table~\ref{tab:three_levels}.

\begin{table}[h!]
\centering
\renewcommand{\arraystretch}{1}
\begin{tabular}{|| c || c | c | c  ||}  
\hline
 logical & $\neg$  & $\land$  &  $\lor$  \\  \hline
 set theoretical & $[]^C$    & $\cap$  &  $\cup$  \\  \hline
algebraic &  $-$    & $*$   &  $+$  \\  \hline
\hline
\end{tabular}
\renewcommand{\arraystretch}{1}
\caption{The purported  interrelationships between three basic operations
at the logic, set, and algebraic levels. .}\label{tab:three_levels}
\end{table}

\section{Multisets in Pattern Recognition}  \label{sec:pattrec}

The possibility to use msets to represent any type of density paves the way to interesting applications
in pattern recognition and deep learning (e.g.~\cite{DudaHart,Koutrombas,Hinton,Schmidhuber,CostaDeep}).  
In this section we illustrate how msets and the Jaccard
index can be readily applied in order to quantify the similarity between two (or more) clusters
represented by respective density functions.  

\begin{figure}[h!]  
\begin{center}
   \includegraphics[width=0.7\linewidth]{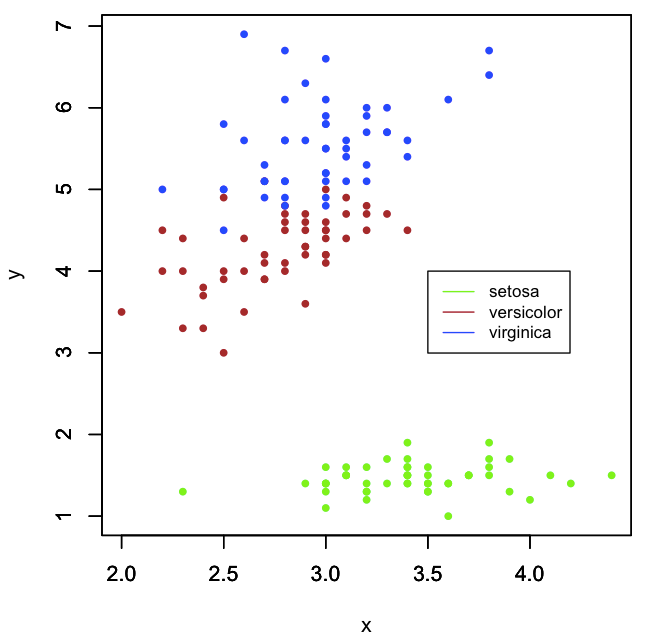}  
    \caption{A scatterplot representing the distribution of three types of iris flowers
    represented by two respective features $x$ and $y$.}
    \label{fig:scatt}
    \end{center}
\end{figure}

Let's consider the three sets of points in the scatterplot shown in Figure~\ref{fig:scatt}, which
corresponds to the three species of iris flower in the frequently adopted iris dataset.  Only two
out of their 4 features have been chosen in the following example for simplicity's sake.

\vspace{0.5cm}

\begin{figure}[h!]  
\begin{center}
   \includegraphics[width=1\linewidth]{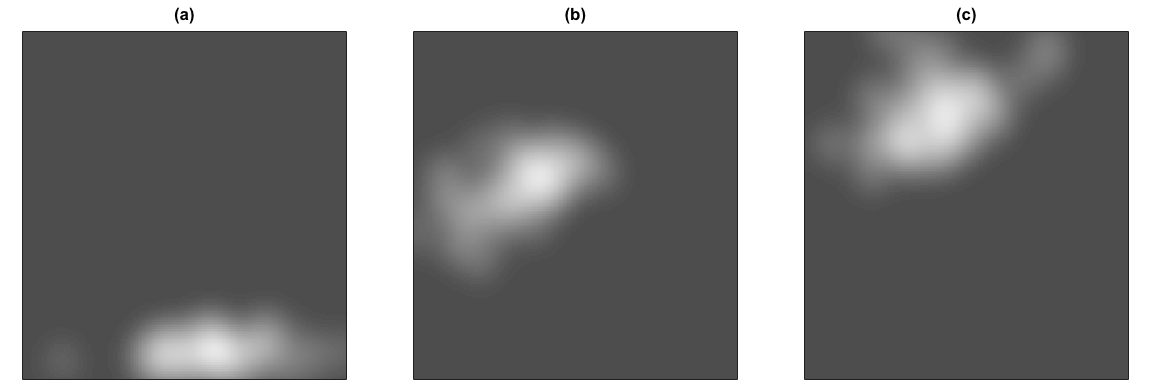}  
    \caption{The three density scalar fields obtained by gaussian kernel expansion
    of each of the three categories.}
    \label{fig:dens}
    \end{center}
\end{figure}

The density obtained from the respective discrete samples through gaussian
kernel expansion are shown in Figure~\ref{fig:dens}.

\vspace{0.5cm}

The obtained multiset Jaccard index for each pairwise combination of categories
are presented in Table~\ref{tab:iris}.

\begin{table}[h!]
\centering
\renewcommand{\arraystretch}{1}
\begin{tabular}{|| c || c | c | c  ||}  
\hline
 & setosa & versicolor & virginica \\
\hline  \hline
setosa & 1 & 2.6e-5 & 0 \\  \hline
versicolor & 2.6e-5 & 1 & 0.145 \\  \hline
viginica & 0 & 0.145 & 1 \\   \hline
\hline
\end{tabular}
\renewcommand{\arraystretch}{1}
\caption{The Jaccard indices obtained for pairwise combinations between the
three iris species.  The resolution has been limited to 6 digits.}\label{tab:iris}
\end{table}

The obtained results are fully compatible with the interrelationships between the
three densities, or clusters, in Figure~\ref{fig:scatt}.  In addition, the threewise
Jaccard index from Equation~\ref{eq:Jac3} result nearly null, indicating a really
small chance that the three densities correspond to the same cluster.

\section{Concluding Remarks}

The fascinating subject of multisets has been presented in a hopefully introductory manner.

We started with a brief review of traditional sets and their properties, which was followed by
a progressive presentation of multisets and their characteristics.  The possibility of obtaining
multiset generalizations capable of dealing with functions, scalar fields, and densities, was
then described and illustrated.

In addition to introducing several of the basic mset concepts, the present work also proposed
how the universe mset can be defined in a robust manner by allowing the subtraction of
msets to take negative values.  This paved the way for recovering several properties analogous
to traditional sets involving the complement operation, including the De Morgan theorem.

The extension of msets to mfunctions, namely traditional real functions represented as
msets,  was also proposed, paving the way to defining mfunctionals, of which the
Jaccard index for mfunctions is one example, including operation analogous to the
traditional inner product, but involving set operations, which paved the way to defining
respective transformations analogous to the Fourier transform, though devoid of orthogonality.

Having defined an operation analogous to the inner product immediately enabled us to
propose binary operators that are mset and mfunction counterparts of the traditional 
convolution between two functions, paving the way to achieving a field that can be
called \emph{integrated signal processing} that is characterized by the incorporation of
mset and mfunction counterparts to most of the concepts and operators adopted in 
signal processing.

The extension of the Jaccard index, which is intrinsically related to set theory, to msets
and also to allow the consideration of more than two sets were also presented, which paves
the way to employing these combined concepts for the characterization of relationships
between clusters in feature spaces, a problem that is common to both the areas of 
pattern recognition and deep learning.

The presented concepts and methods can lead to several interesting applications,
also motivating further integrations between the structures and properties between the
domains of set theory, propositional logic, and analysis.

\vspace{0.7cm}
\emph{Acknowledgments.}

Luciano da F. Costa
thanks CNPq (grant no.~307085/2018-0) and FAPESP (grant 15/22308-2).  
\vspace{1cm}

\bibliography{mybib}

\begin{thebibliography}{10}

\bibitem{Hein}
J.~Hein.
\newblock {\em Discrete Mathematics}.
\newblock Jones \& Bartlett Pub., 2003.

\bibitem{Knuth}
D.~E. Knuth.
\newblock {\em The Art of Computing}.
\newblock Addison Wesley, 1998.

\bibitem{Blizard}
W.~D. Blizard.
\newblock Multiset theory.
\newblock {\em Notre Dame Journal of Formal Logic}, 30:36---66, 1989.

\bibitem{Blizard2}
W.~D. Blizard.
\newblock The development of multiset theory.
\newblock {\em Modern Logic}, 4:319?352, 1991.

\bibitem{Thangavelu}
P.~M. Mahalakshmi and P.~Thangavelu.
\newblock Properties of multisets.
\newblock {\em International Journal of Innovative Technology and Exploring
  Engineering}, 8:1--4, 2019.

\bibitem{Singh}
D.~Singh, M.~Ibrahim, T.~Yohana, and J.~N. Singh.
\newblock Complementation in multiset theory.
\newblock {\em International Mathematical Forum}, 38:1877--1884, 2011.

\bibitem{DudaHart}
R.~O. Duda, P.~E. Hart, and D.~G. Stork.
\newblock {\em Pattern Classification}.
\newblock Wiley Interscience, 2000.

\bibitem{Koutrombas}
K.~Koutrombas and S.~Theodoridis.
\newblock {\em Pattern Recognition}.
\newblock Academic Press, 2008.

\bibitem{Hinton}
G.~E. Hinton.
\newblock Training products of experts by mini-mizing contrastive divergence.
\newblock {\em Neural computation}, 14(8):1771--1800, 2002.

\bibitem{Schmidhuber}
J.~Schmidhuber.
\newblock Deep learning in neural networks:an overview.
\newblock {\em Neural networks}, 61:85--117, 2015.

\bibitem{CostaDeep}
H.~F. de~Arruda, A.~Benatti, C.~H. Comin, and L.~da~F.~Costa.
\newblock Learning deep learning.
\newblock Researchgate, 2019.
\newblock
  \url{https://www.researchgate.net/publication/335798012_Learning_Deep_Learning_CDT-15}.
  [Online; accessed 22-Dec-2019.].

\bibitem{CostaModeling}
L.~da~F.~Costa.
\newblock Modeling: The human approach to science.
\newblock Researchgate, 2019.
\newblock
  \url{https://www.researchgate.net/publication/333389500_Modeling_The_Human_Approach_to_Science_CDT-8}.
  [Online; accessed 1-Oct-2020.].

\bibitem{CostaAmple}
L.~da~F.~Costa.
\newblock An ample approach to modeling.
\newblock Researchgate, 2019.
\newblock
  \url{https://www.researchgate.net/publication/355056285_An_Ample_Approach_to_Data_and_Modeling}.
  [Online; accessed 10-Oct-2021.].

\bibitem{Jaccard}
P.~Jaccard.
\newblock {\'E}tude comparative de la distribution florale dans une portion des
  alpes et des jura.
\newblock {\em Bulletin de la Soci\'et\'e vaudoise des sciences naturelles},
  37:547--549, 1901.

\bibitem{Jac:wiki}
Wikipedia.
\newblock Jaccard index.
\newblock \url{https://en.wikipedia.org/wiki/Jaccard_index}. [Online; accessed
  10-Oct-2021].

\bibitem{CostaJaccard}
L.~da~F. Costa.
\newblock Further generalizations of the jaccard index.
\newblock
  \url{https://www.researchgate.net/publication/355381945_Further_Generalizations_of_the_Jaccard_Index},
  2021.
\newblock [Online; accessed 21-Aug-2021].

\bibitem{Harmuth}
H.~F. Harmuth.
\newblock Applications of walsh functions in communications.
\newblock {\em IEEE Spectrum}, 6:82--91, 1969.

\bibitem{Tzafestas}
S.~G. Tzafestas.
\newblock {\em Walsh Functions in Signal and Systems Analysis and Design}.
\newblock Van Nostrand Reinhold, New York, 1985.

\end{thebibliography}
\bibliographystyle{unsrt}

\end{document}